\documentclass[11pt,reqno]{amsart}
\usepackage{amsmath,amssymb,graphicx}
\newtheorem{thm}{Theorem}[section]
\newtheorem{cor}[thm]{Corollary}
\newtheorem{lem}[thm]{Lemma}
\newtheorem{prop}[thm]{Proposition}
\theoremstyle{definition}
\newtheorem{defn}[thm]{Definition}
\theoremstyle{remark}
\newtheorem{rem}[thm]{Remark}
\theoremstyle{remark}
\newtheorem{ex}[thm]{Example}
\numberwithin{equation}{section}

\begin{document}

\title{Non stopping times and stopping theorems}%
\author{Ashkan Nikeghbali}
\address{ETHZ \\ Departement Mathematik, R\"{a}mistrasse 101, HG G16\\ Z\"{u}rich 8092, Switzerland}
 \email{ashkan.nikeghbali@math.ethz.ch}
\date{\today}
%\dedicatory{}%
%\commby{}%
% ----------------------------------------------------------------
\begin{abstract}
 Given a random time, we give some characterizations of the set of martingales for which the stopping
theorems still hold. We also investigate how the stopping theorems
are modified when we consider arbitrary random times. To this end,
we introduce some families of martingales with remarkable
properties.
\end{abstract}
\maketitle
% ----------------------------------------------------------------

\section{Introduction}
The role of stopping times in martingale theory is fundamental. In
particular, there are myriads of applications of Doob's optional
stopping theorem:

\begin{itemize}
\item If $\left( M_{t}\right) $ is a uniformly integrable martingale, and $T$%
\ is a stopping time (both with respect to the filtration $\left( \mathcal{F}%
_{t}\right) $ which is assumed to satisfy the usual hypotheses under
a given
probability space $\left( \Omega ,\mathcal{F},\mathbb{P}\right) $), then:%
\begin{equation}
\mathbb{E}\left[ M_{T}\right] =\mathbb{E}\left[ M_{\infty }\right] =\mathbb{E%
}\left[ M_{0}\right]  \label{optstopthm}
\end{equation}%
and, in fact:%
\begin{equation}
\mathbb{E}\left[ M_{\infty }\mid \mathcal{F}_{T}\right] =M_{T}
\label{KM}
\end{equation}
\end{itemize}

In this paper, we would like to discuss in some depth the following
question, which arises very naturally:

\begin{itemize}
\item What happens to (\ref{optstopthm}) and (\ref{KM}) when $T$\ is
replaced by a random time $\rho$, and $\mathcal{F}_{T}$\ by $\mathcal{F}%
_{\rho}=\sigma \left\{ H_{\rho};\text{ }H\text{ is }\left(
\mathcal{F}_{t}\right) \mathcal{\ }\text{optional}\right\} $?
\end{itemize}

Some partial answers to this general question have been given by D.
Williams \cite{williams} on one hand, and Knight-Maisonneuve
\cite{Knight-Maisonneuve} on the other hand:

\begin{enumerate}
\item There exist \textquotedblleft non-stopping\textquotedblright\ times $%
\rho$, which we have called pseudo-stopping times in
\cite{AshkanYor} such that, for every bounded martingale $M$,
\begin{equation}
\mathbb{E}\left[ M_{\rho}\right] =\mathbb{E}\left[ M_{0}\right]
\label{pdtacond}
\end{equation}

\item If for every bounded martingale $M$, one has:%
\begin{equation}
\mathbb{E}\left[ M_{\infty }\mid \mathcal{F}_{\rho}\right]
=M_{\rho}, \label{KMcond}
\end{equation}%
then $\rho$\ is a stopping time (\cite{Knight-Maisonneuve}).
\end{enumerate}

Besides the fact that it is mathematically an interesting question
to understand how and why the usual results fail to hold when
stopping times are replaced with arbitrary random times, it should
be noticed that random times that are not stopping times play a key
role in various contexts, such as in the modeling of default times
in mathematical finance (see \cite{elliotjeanbyor}), in Markov
Processes theory (see \cite{delmaismey}), in the characterization of
the set of zeros of continuous martingales (\cite{azemayorzero}), in
path decomposition of some diffusions (see \cite{jeulin} or
\cite{ashyordoob}), in the study of Strong Brownian Filtrations (see
\cite{beksy}), etc.

The most studied family of random times, after stopping times, are
ends of optional sets, also named honest times (such the last zero
of the standard Brownian Motion before a fixed time). A very
powerful, but not so well known, technique to study such random
times, is the progressive expansions or enlargements of filtrations.
The theory of progressive enlargements of filtrations was introduced
independently by Barlow (\cite{barlow}) and Yor (\cite{gros}), and
further developed by Jeulin and Yor
(\cite{jeulinyor,yorjeulin,jeulin,zurich}). The reader can find many
applications of this theory in the cited references and in
\cite{AshkanYor} and \cite{ashyordoob}. The concept of dual
projections also play an important role in the study of arbitrary
random times (see \cite{dellachmeyer} or \cite{rogerswilliams}).

The main idea in the progressive enlargements setting is to consider
the larger filtration $\left(\mathcal{F}_{t}^{\rho}\right)$, which
is the smallest right continuous filtration which contains
$\left(\mathcal{F}_{t}\right)$ and which makes $\rho$ a stopping
time, and then to see how martingales of the smaller filtration are
changed when considered as processes of the larger one. In
\cite{azemajeulinknightyor}, the authors used these ideas to give a
solution to equation (\ref{KMcond}) in a Brownian setting, using a
predictable representation property for martingales in the larger
filtration $\left(\mathcal{F}_{t}^{\rho}\right)$. In this paper, we
shall solve equation (\ref{pdtacond}) for arbitrary random times and
equation (\ref{KMcond}) for honest times. In this latter case, we
propose two different approaches and our characterizations ((Theorem
\ref{resolutionune} and Proposition \ref{resolutiondeux})) of the
set of martingales which satisfy (\ref{KMcond}) will be different
(but not necessarily more handy) from the one in
\cite{azemajeulinknightyor}, in that our solution is based only on
quantities relative to the filtration
$\left(\mathcal{F}_{t}\right)$, which moreover is not assumed to be
a Brownian filtration. More precisely, the organization of the paper
is as follows:

\noindent In Section 2, we recall some basic facts about progressive
enlargements of filtrations and arbitrary random times.

\noindent In Section 3, we solve equation (\ref{pdtacond}) for
arbitrary random times, using elementary properties of dual
projections and Laguerre polynomials .

\noindent In Section 4, we provide two different approaches to solve
(\ref{KMcond}) for honest times (Theorem \ref{resolutionune} and
Proposition \ref{resolutiondeux}). In particular, we shall see how
to obtain a large class of solutions to (\ref{KMcond}), by
considering martingales which vanish at $L$. We illustrate these
facts in the celebrated special case when $L$ is the last time
before a fixed (or a stopping) time when a standard Brownian Motion
vanishes.

\noindent In Section 5, we introduce a family of test martingales,
with interesting and universal properties (in a sense that will be
clear), to understand how the equalities (\ref{pdtacond}) and
(\ref{KMcond}) may fail to hold for honest times.
\section{Basic facts about progressive enlargements of filtrations
and random times} In this Section, we recall some results (which may
not be so well known) that we shall use in this paper and fix the
notations once and for all. Throughout this article, we assume for
simplicity that $\rho$ is a random time such that
$\mathbb{P}\left[\rho=0\right]=\mathbb{P}\left[\rho=\infty\right]=0$.

Let $\left( \Omega ,\mathcal{F},\left( \mathcal{F}_{t}\right) _{t\geq 0},%
\mathbb{P}\right) $ a filtered probability space, and $\rho :$
$\left(
\Omega ,\mathcal{F}\right) \rightarrow \left( \mathbb{R}_{+},\mathcal{B}%
\left( \mathbb{R}_{+}\right) \right) $ be a random time. We enlarge
the initial filtration $\left( \mathcal{F}_{t}\right) $\ with the
process $\left( \rho \wedge t\right) _{t\geq 0}$, so that the new
enlarged filtration $\left( \mathcal{F}_{t}^{\rho }\right) _{t\geq
0}$\ is the
smallest filtration containing $\left( \mathcal{F}_{t}\right) $\ and making $%
\rho $\ a stopping time. A few processes will play a crucial role in
our discussion:

\begin{itemize}
\item the $\left( \mathcal{F}_{t}\right) $-supermartingale
\begin{equation}
Z_{t}^{\rho }=\mathbb{P}\left[ \rho >t\mid \mathcal{F}_{t}\right]
\label{surmart}
\end{equation}%
chosen to be c\`{a}dl\`{a}g, associated to $\rho $\ by Az\'{e}ma (see \cite%
{jeulin} for detailed references);

\item the $\left( \mathcal{F}_{t}\right) $ dual optional
projection of the process $1_{\left\{ \rho \leq t\right\} }$,
denoted by $A_{t}^{\rho }$;

\item the c\`{a}dl\`{a}g martingale
\begin{equation*}
\mu _{t}^{\rho }=\mathbb{E}\left[ A_{\infty }^{\rho }\mid \mathcal{F}_{t}%
\right] =A_{t}^{\rho }+Z_{t}^{\rho }
\end{equation*}%
which is in BMO($\mathcal{F}_{t}$) (see \cite{delmaismey} or
\cite{zurich}).
\end{itemize}
Every $\left( \mathcal{F}_{t}\right) $ local
martingale $\left( M_{t}\right) $, stopped at $\rho $, is a $\left( \mathcal{%
F}_{t}^{\rho }\right) $\ semimartingale, with canonical
decomposition:
\begin{equation}
M_{t\wedge \rho }=\widetilde{M}_{t}+\int_{0}^{t\wedge \rho
}\dfrac{d\langle M,\mu ^{\rho }\rangle_{s}}{Z_{s-}^{\rho }}
\label{decocanonique}
\end{equation}%
where $\left( \widetilde{M}_{t}\right) $\ is an $\left( \mathcal{F}%
_{t}^{\rho }\right) $-local martingale.

The most interesting case in the theory of progressive enlargements
of filtrations is when $\rho$ is an honest time; we will always
denote honest times by $L$ instead of $\rho$. Indeed, if $L$ is an
honest time, then every $\left( \mathcal{F}_{t}\right) $ local
martingale $\left( M_{t}\right) $, is an $\left( \mathcal{%
F}_{t}^{\rho }\right) $\ semimartingale, with canonical
decomposition:
\begin{equation}
M_{t}=\widetilde{M}_{t}+\int_{0}^{t\wedge L }\frac{d\langle M,\mu
^{L }\rangle_{s}}{Z_{s-}^{L }}-\int_{L}^{t}\dfrac{d\langle M,\mu ^{L
}\rangle_{s}}{1-Z_{s-}^{L}}. \label{decocanoniquehonete}
\end{equation}We shall often need to make one (or sometimes both) of
the following assumptions:
\begin{itemize}
\item Assumption $\mathbf{(C)}$: all $\left( \mathcal{F}_{t}\right) $-martingales are \underline{c}ontinuous (e.g:
 the Brownian filtration).
 \item Assumption $\mathbf{(A)}$: the random time $\rho $ \underline{a}voids every $\left( \mathcal{F}_{t}\right) $%
-stopping time $T$, i.e. $\mathbb{P}\left[ L =T\right] =0$.
\end{itemize}When we refer to assumptions $\mathbf{(CA)}$, this will
mean that both the conditions $\mathbf{(C)}$ and $\mathbf{(A)}$
hold. Under conditions $\mathbf{(C)}$ or $\mathbf{(A)}$,
$A_{t}^{\rho }$ is also the dual predictable projection of
$1_{\left\{ \rho \leq t\right\} }$; moreover under $\mathbf{(A)}$,
$A_{t}^{\rho }$ is continuous.

Now, we give the definitions of some sigma fields associated with
arbitrary random times, following Chung and Doob (\cite{chungdoob}):
\begin{defn}
Three classical $\sigma $-fields associated with a filtration
$\left( \mathcal{F}_{t}\right) \ $and any random time $\rho $ are:
\begin{equation*}
\left\{
\begin{array}{lcl}
\mathcal{F}_{\rho +} & = & \sigma \left\{ z_{\rho },\text{ }\left(
z_{t}\right) \text{ any }\left( \mathcal{F}_{t}\right) \
\text{progressively
measurable process}\right\} ; \\
\mathcal{F}_{\rho } & = & \sigma \left\{ z_{\rho },\text{ }\left(
z_{t}\right) \text{ any }\left( \mathcal{F}_{t}\right) \
\text{optional
process}\right\} ; \\
\mathcal{F}_{\rho -} & = & \sigma \left\{ z_{\rho },\text{ }\left(
z_{t}\right) \text{ any }\left( \mathcal{F}_{t}\right) \
\text{predictable
process}\right\} ;%
\end{array}%
\right.
\end{equation*}
\end{defn}
Under condition $\mathbf{(A)}$, we have: $\mathcal{F}_{\rho
}=\mathcal{F}_{\rho -}$.

We conclude this section by giving two results (due to Az\'{e}ma
\cite{azema}) which will play an important role in the next
sections. The reader can also refer to the book \cite{delmaismey}
for a very nice introduction to the results of Az\'{e}ma and the
theory of progressive enlargements of filtrations.
\begin{lem}[Az\'{e}ma \cite{azema}]\label{azemgeneral}
Let $L$ be an honest time; then under $(\mathbf{A})$,
$A_{\infty}^{L}$ follows the exponential law with parameter $1$ and
the measure $dA_{t}^{L}$ is carried by the set
$\left\{t:\;Z_{t}^{L}=1\right\}$. Moreover, $A^{L}$ does not
increase after $L$, i.e. $A_{L}^{L}=A_{\infty}^{L}$.
\end{lem}
\begin{lem}[Az\'{e}ma \cite{azema}]\label{azemalemmdeux}
Let $L$ be an honest time and assume $(\mathbf{A})$ holds. Then,
$$L=\sup\left\{t:\;1-Z_{t}^{L}=0\right\}.$$In particular, $1-Z_{L}^{L}=0$.
\end{lem}

\section{A resolution of the equation $\mathbb{E}\left( M_{\rho}\right) =\mathbb{E}\left(
M_{0}\right)$}
 We wish to solve
equation (\ref{pdtacond}), where $\rho$ is given and the unknown are
all bounded $\left( \mathcal{F}_{t}\right)$ martingales for which
(\ref{pdtacond}) hold. We consider the class of bounded martingales
because we want to make sure that $\mathbb{E}\left( M_{\rho}\right)$
exists; in fact, we can look for solutions to equation
(\ref{pdtacond}) in the space of $\mathcal{H}^{1}$ martingales (see
\cite{kazamaki} or \cite{meyerdual}). We recall that the space
$\mathcal{H}^{1}$ is the Banach space of (c\`{a}dl\`{a}g) $\left(
\mathcal{F}_{t}\right) $ martingales $\left( M_{t}\right) $ such
that
\begin{equation*}
\left\Vert M\right\Vert _{\mathcal{H}^{1}}=\mathbb{E}\left[
\sup_{t\geq 0}\left\vert M_{t}\right\vert \right] <\infty .
\end{equation*}
\begin{defn}
We call $\mathcal{S}_{1}$ the set of solutions of equation
(\ref{pdtacond}),
i.e.$$\mathcal{S}_{1}\equiv\left\{M\in\mathcal{H}^{1}:\;\mathbb{E}\left(
M_{\rho}\right)=\mathbb{E}\left( M_{\infty}\right)\right\}.$$
\end{defn}
\begin{thm}\label{caracter1}
The map
$$T\left(M\right)=\mathbb{E}\left[\langle M,\mu^{\rho}\rangle_{\infty}\right],$$defines
a continuous linear form on the Banach space $\mathcal{H}^{1}$, and
we have the following characterizations for $\mathcal{S}_{1}$:
\begin{enumerate}
\item $$\mathcal{S}_{1}=\ker T,$$or in other words,
\begin{equation} \label{condnoyau}
\mathcal{S}_{1}=\left\{M\in\mathcal{H}^{1}:\;\mathbb{E}\left[\langle
M,\mu^{\rho}\rangle_{\infty}\right]=0\right\}.
\end{equation}
\item \begin{equation} \label{dualhbmo}
\mathcal{S}_{1}=\left\{M\in\mathcal{H}^{1}:\;\mathbb{E}\left[M_{\infty}\left(\mu^{\rho}_{\infty}-1\right)\right]=0\right\}.
\end{equation}
\item \begin{equation} \label{projdualprev}
\mathcal{S}_{1}=\left\{M\in\mathcal{H}^{1}:\;\mathbb{E}\left[M_{\infty}\left(A^{\rho}_{\infty}-1\right)\right]=0\right\}.
\end{equation}
\end{enumerate}Consequently, $\mathcal{S}_{1}$ is a closed linear subspace of $\mathcal{H}^{1}$.
\end{thm}
\begin{proof}
The fact that $T$ defines a linear form is a consequence of the well
known duality between $\mathcal{H}^{1}$ and $BMO$ (see
\cite{kazamaki} or \cite{meyerdual} for details and references).
Now, let $\rho$ be a random time and let $M$ be a martingale in
$\mathcal{H}^{1}$. We have (see \cite{dellachmeyer} or
\cite{Ashkanassay} where dual projections and their properties
applications are discussed):
\begin{equation}\label{aa}
    \mathbb{E}\left[M_{\rho}\right]=\mathbb{E}\left[\int_{0}^{\infty}M_{s}dA_{s}^{\rho}\right]=\mathbb{E}\left[M_{\infty}A_{\infty}^{\rho}\right].
\end{equation}Hence, $$\mathbb{E}\left[M_{\rho}\right]=\mathbb{E}\left[M_{\infty}\right]\Leftrightarrow \mathbb{E}\left[M_{\infty}\left(A^{\rho}_{\infty}-1\right)\right]=0,$$
and this establishes $(3)$.

But as $\mathbb{P}\left[\rho=\infty\right]=0$,
$Z_{\infty}^{\rho}=0$, and $\mu_{\infty}^{\rho}=A_{\infty}^{\rho}$.
We thus have:
\begin{equation*}
\mathbb{E}\left[M_{\rho}\right]=\mathbb{E}\left[M_{\infty}\mu_{\infty}^{\rho}\right]=\mathbb{E}\left[M_{\infty}\right]+\mathbb{E}\left[\langle
M,\mu^{\rho}\rangle_{\infty}\right],
\end{equation*}and  $(1)$ and $(2)$ follow easily.
\end{proof}
\begin{rem}
As will be shown later, one must not confuse (\ref{condnoyau}) with
the stronger condition $\langle M,\mu^{\rho}\rangle_{t}=0$ for every
$t$.
\end{rem}
Theorem \ref{caracter1} shows that the set of solutions of equation
(\ref{pdtacond}) is a linear space of codimension $1$ in
$\mathcal{H}^{1}$ if the linear form $T$ is not null. The case when
this form is null corresponds to the remarkable class of random
times called pseudo-stopping times, defined and studied in
\cite{AshkanYor}.
\begin{prop}[\cite{AshkanYor}]
The following are equivalent:
\begin{enumerate}
\item (\ref{pdtacond}) holds for every martingale in
$\mathcal{H}^{1}$;
\item $A_{\infty}^{\rho}=\mu_{\infty}^{\rho}=1\;\mathrm{a.s.}$;
\item If $\left(M_{t}\right)$ is an $\left(\mathcal{F}_{t}\right)$
local martingale, then $\left(M_{t\wedge \rho}\right)$ is an
$\left(\mathcal{F}_{t}^{\rho}\right)$ local martingale.
\end{enumerate}
\end{prop}
We can also give the following elementary but useful corollary of
Theorem \ref{caracter1}:
\begin{cor}\label{relaorthogonalite}
Let $M$ be an $L^{2}$ bounded martingale such that $M_{\infty}\in
\left(L^{2}\left(\sigma\left(A_{\infty}^{\rho}\right)\right)\right)^{\perp}$,
the orthogonal of $
L^{2}\left(\sigma\left(A_{\infty}^{\rho}\right)\right)$. Then
$M\in\mathcal{S}_{1}$.
\end{cor}

\bigskip
Now, one may want to find some $L^{2}$ bounded martingales such that
$M_{\infty}\in
L^{2}\left(\sigma\left(A_{\infty}^{\rho}\right)\right)$. This can be
done with the help of orthogonal polynomials if one knows the law of
$A_{\infty}^{\rho}$. We shall now illustrate this  with the
important case of honest times,  giving a complete description of
$\mathcal{S}_{1}$ in terms of Laguerre polynomials.
\bigskip

We first introduce some basic facts about Laguerre polynomials
(\cite{AARoy}). Let us consider the Hilbert space $L^{2}\left( \exp
\left( -x\right) dx\right) $. The Laguerre Polynomials,
$\widetilde{\mathcal{L}}_{n}\left( x\right) $ are the orthogonal
polynomials associated with the measure $\exp \left(
-x\right) dx$. They are given by the formula:%
\begin{eqnarray*}
\widetilde{\mathcal{L}}_{n}\left( x\right) &=&\sum_{k=0}^{n}\left( -1\right) ^{k}\frac{%
\left( n!\right) ^{2}}{\left( k!\right) ^{2}\left( n-k\right) !}x^{k} \\
&=&\exp \left( x\right) \frac{d^{n}}{dx^{n}}\left( x^{n}\exp \left(
-x\right) \right)
\end{eqnarray*}%
They satisfy the orthogonality relation:%
\begin{equation*}
\int_{0}^{\infty }dx\exp \left( -x\right)
\widetilde{\mathcal{L}}_{m}\left( x\right)
\widetilde{\mathcal{L}}_{n}\left( x\right) =\left( n!\right)
^{2}\delta _{m,n}
\end{equation*}%
We can normalize and take:%
\begin{equation*}
\mathcal{L}_{n}\left( x\right)
=\frac{\widetilde{\mathcal{L}}_{n}\left( x\right) }{n!}
\end{equation*}%
so that the family $\left( \mathcal{L}_{n}\left( x\right) \right) $
is an orthonormal basis in $L^{2}\left( \exp \left( -x\right)
dx\right) $. For example,
\begin{eqnarray*}
\mathcal{L}_{0}\left( x\right) &=&1 \\
\mathcal{L}_{1}\left( x\right) &=&1-x \\
\mathcal{L}_{2}\left( x\right) &=&\frac{1}{2}\left( 2-4x+x^{2}\right) \\
\mathcal{L}_{3}\left( x\right) &=&\frac{1}{6}\left( 6-18x+9x^{2}-x^{3}\right) \\
\mathcal{L}_{4}\left( x\right) &=&\frac{1}{24}\left(
24-96x+72x^{2}-16x^{3}+x^{4}\right).
\end{eqnarray*}
\begin{thm}
Let $L$ be an honest time and assume condition $(\mathbf{A})$ holds.
Let $M$ be an $L^{2}$ bounded martingale. Then the following are
equivalent:
\begin{enumerate}
\item $M\in \mathcal{S}_{1}$;
\item $M_{\infty}$ may be represented as:
\begin{equation}\label{lag}
    M_{\infty}=X+\varphi\left(A_{\infty}\right),
\end{equation}where $X\in\left(L^{2}\left(\sigma\left(A_{\infty}\right)\right)\right)^{\perp}$
and where $\varphi\in
L^{2}\left(\sigma\left(A_{\infty}\right)\right)$ admits the
following representation:
\begin{equation}\label{represdephi}
\varphi\left(A_{\infty}\right)=\alpha_{0}+\sum_{n=2}^{\infty }\alpha
_{n}\mathcal{L}_{n}\left( A_{\infty}\right),
\end{equation}
 with $\alpha\in \mathbb{R}$ and $\left(
\alpha _{n}\right)$ such that $\sum \alpha _{n}^{2}<\infty $, i.e:
in the development of $\varphi$, the coefficient of $\mathcal{L}_{1}
$ is $\alpha _{1}=0$.
\end{enumerate}
\end{thm}
\begin{proof}
We note that:
$$M_{\infty}=X+\mathbb{E}\left[M_{\infty}|A_{\infty}\right]\equiv X+\varphi\left(A_{\infty}\right),$$ with
$X\in\left(L^{2}\left(\sigma\left(A_{\infty}\right)\right)\right)^{\perp}$
and $\varphi\left(A_{\infty}\right) \in
L^{2}\left(\sigma\left(A_{\infty}\right)\right)$. Now, from
(\ref{projdualprev}), $M\in \mathcal{S}_{1}$ if and only if
$$\mathbb{E}\left[\varphi\left(A_{\infty}\right)\left(A_{\infty}-1\right)\right]=0,$$or
equivalently
\begin{equation}\label{lignedessus}
\mathbb{E}\left[\varphi\left(A_{\infty}\right)\mathcal{L}_{1}\left(A_{\infty}\right)\right]=0.
\end{equation}
Since the family $\left(\mathcal{L}_{n}\right)_{n\geq 0}$ is total
in $L^{2}\left(\exp\left(-x\right)dx\right)$, we can represent
$\varphi\left(x\right)$ as:
$$\varphi\left(x\right)
=\sum_{n=0}^{\infty}\alpha_{n}\mathcal{L}_{n}\left(x\right),$$with
$\sum \alpha _{n}^{2}<\infty $. Now putting the series expansion of
$\varphi$ in (\ref{lignedessus}) and using the fact that the family
$\left(L_{n}\right)_{n\geq 0}$ is orthogonal gives the desired
result.
\end{proof}
\begin{ex}
Let the filtration $\left( \mathcal{F}_{t}\right) $\ be generated by
a one dimensional Brownian motion $\left( B_{t}\right) _{t\geq 0}$,
and let
\begin{equation*}
T_{1}=\inf \left\{ t:\text{ }B_{t}=1\right\} ,\text{ and }L =\sup
\left\{ t<T_{1}:\text{ }B_{t}=0\right\}.
\end{equation*}It is well known that
$$Z_{t}=\mathbb{P}\left[L>t\mid\mathcal{F}_{t}\right]=1-B^{+}_{t\wedge
T_{1}}$$An application of Tanaka's formula yields:
$A_{\infty}=\dfrac{1}{2}\ell_{T_{1}}$, where $\left(\ell_{t}\right)$
is the Brownian local time at zero. Now, from the previous theorem,
it is easily seen that any martingale of the form
$M_{t}=\mathbb{E}\left[\mathcal{L}_{n}\left(\ell_{T_{1}}\right)\mid\mathcal{F}_{t}\right];\;n\neq
1$ is in $\mathcal{S}_{1}$.
\end{ex}\bigskip
It is also possible to use the Kunita-Watanabe orthogonal
decompositions for square integrable martingales to give a
description of $\mathcal{S}_{1}$; more precisely:
\begin{prop}
Let $M$ be an $L^{2}$ martingale and let $\rho$ be an arbitrary
random time. Then the following are equivalent:
\begin{enumerate}
\item $M\in \mathcal{S}_{1}$;
\item $\left(M_{t}\right)$ decomposes as:
\begin{equation}\label{kunitwatanbe}
    M_{t}=N_{t}+\int_{0}^{t}k_{s}d\mu_{s}^{\rho},
\end{equation}where $N$ is an $L^{2}$ martingale such that
$$\langle N,\mu^{\rho}\rangle_{t}=0,\;\forall t\geq 0$$ and $k$ is a predictable
process such that:
$$\mathbb{E}\left[\int_{0}^{\infty}k_{s}^{2}d\langle\mu^{\rho}\rangle_{s}\right]<\infty;\;\mathbb{E}\left[\int_{0}^{\infty}k_{s}d\langle \mu^{\rho}\rangle_{s}\right]=0.$$
\end{enumerate}
\end{prop}
\begin{proof}
From the Kunita-Watanbe decomposition (\cite{Kunita-Watanabe}), any
$L^{2}$ martingale $M$ can be decomposed as:
$M_{t}=N_{t}+\int_{0}^{t}k_{s}d\mu_{s}^{\rho}$, where $\langle
N,\mu^{\rho}\rangle_{t}=0,\;\forall t\geq 0,$ and where $k$ is a
predictable process such that
$\mathbb{E}\left[\int_{0}^{\infty}k_{s}^{2}d\langle\mu^{\rho}\rangle_{s}\right]<\infty$.
Now, from (\ref{condnoyau}), it follows that $M\in \mathcal{S}_{1}$
if and only if
$$\mathbb{E}\left[\langle M,\mu^{\rho}\rangle_{\infty}\right]=0,$$
or equivalently
$$\mathbb{E}\left[\int_{0}^{\infty}k_{s}d\langle\mu^{\rho}\rangle_{s}\right]=0,$$which
completes the proof.
\end{proof}
\section{A resolution of the equation $\mathbb{E}\left( M_{\infty }\mid \mathcal{F}_{L}\right)
=M_{L}$}In this Section, we shall try to give explicit solutions to
equation (\ref{KMcond}), when $\rho\equiv L$ is an honest time
satisfying $(\mathbf{A})$.
\begin{defn}
We call $\mathcal{S}_{2}$ the set of solutions to equation
(\ref{KMcond}): $$\mathcal{S}_{2}=\left\{M\in
\mathcal{H}^{1}:\;\mathbb{E}\left[ M_{\infty }\mid
\mathcal{F}_{L}\right] =M_{L}\right\}.$$
\end{defn}
\begin{rem}
$\mathcal{S}_{2}\subset \mathcal{S}_{1}$.
\end{rem}

We recall here that equation (\ref{KMcond}) was solved, in the case
of the Brownian filtration, by Az\'{e}ma, Knight, Jeulin and Yor in
\cite{azemajeulinknightyor}. We propose two other characterizations
of the set of solutions to this equation (Theorem
\ref{resolutionune} and Proposition \ref{resolutiondeux}). From now
on, we assume that $L$ is an honest time satisfying $(\mathbf{A})$.
\subsection{A general solution related to the enlargements formulae}
Recall that under condition $(\mathbf{A})$, $\mathcal{F}_{L
}=\mathcal{F}_{L -}$.
\begin{lem}\label{lemtribus}
$\mathcal{F}_{L }=\mathcal{F}_{L- }^{L}$.
\end{lem}
\begin{proof}
From results of Jeulin (\cite{jeulin}), every
$\left(\mathcal{F}_{t}^{L}\right)$ predictable process $H$ can be
represented as
$$H=J\mathbf{1}_{\left]0,L\right]}+K\mathbf{1}_{\left]L,\infty\right[},$$where
$J$ and $K$ are $\left(\mathcal{F}_{t}\right)$ predictable
processes. Hence, $\mathcal{F}_{L- }^{L}=\mathcal{F}_{L -}$, and
since under $(\mathbf{A})$, $\mathcal{F}_{L }=\mathcal{F}_{L -}$,
the lemma is proved.
\end{proof}
\begin{rem}
Using the representation of optional
$\left(\mathcal{F}_{t}^{L}\right)$ processes, it is possible to show
that $\mathcal{F}_{L+}=\mathcal{F}_{L}^{L}$ (see \cite{jeulin}).
\end{rem}
Now, we state a general necessary and sufficient condition for $M$
to be in $\mathcal{S}_{2}$.
\begin{thm}\label{resolutionune}
Let $M\in \mathcal{H}^{1}$. The following are equivalent:
\begin{enumerate}
\item $M\in\mathcal{S}_{2}$;
\item $$\mathbb{E}\left[\int_{0}^{\infty}\dfrac{d\langle M,\mu\rangle_{s}}{1-Z_{s-}}\mid\mathcal{F}_{L
}\right]=\int_{0}^{L}\dfrac{d\langle
M,\mu\rangle_{s}}{1-Z_{s-}},$$or
equivalently:$$\mathbb{E}\left[\int_{L}^{\infty}\dfrac{d\langle
M,\mu\rangle_{s}}{1-Z_{s-}}\mid\mathcal{F}_{L }\right]=0.$$
\end{enumerate}
\end{thm}
\begin{proof}
Let $M\in\mathcal{H}^{1}$; then from (\ref{decocanoniquehonete}),
there exists an $\mathcal{F}_{t}^{L}$ martingale $\widetilde{M}$
such that:
\begin{equation*}
M_{t}=\widetilde{M}_{t}+\int_{0}^{t\wedge L }\frac{d\langle M,\mu
\rangle_{s}}{Z_{s-}}-\int_{L}^{t}\dfrac{d\langle
M,\mu\rangle_{s}}{1-Z_{s-}}.
\end{equation*}We deduce from this decomposition formula
that:
\begin{equation}\label{e1}
M_{L}=\widetilde{M}_{L}+\int_{0}^{L}\frac{d\langle M,\mu
\rangle_{s}}{Z_{s-}},
\end{equation}and
\begin{equation}\label{e2}
\mathbb{E}\left[ M_{\infty }\mid
\mathcal{F}_{L}\right]=\mathbb{E}\left[ \widetilde{M}_{\infty }\mid
\mathcal{F}_{L}\right]+\int_{0}^{L}\frac{d\langle M,\mu
\rangle_{s}}{Z_{s-}}+\mathbb{E}\left[\int_{L}^{\infty}\dfrac{d\langle
M,\mu\rangle_{s}}{1-Z_{s-}}\mid\mathcal{F}_{L }\right].
\end{equation}Now, from Lemma \ref{lemtribus}, $$\mathbb{E}\left[ \widetilde{M}_{\infty }\mid
\mathcal{F}_{L}\right]=\mathbb{E}\left[ \widetilde{M}_{\infty }\mid
\mathcal{F}_{L-}^{L}\right]=\mathbb{E}\left[
\mathbb{E}\left[\widetilde{M}_{\infty
}\mid\mathcal{F}_{L}^{L}\right]\mid \mathcal{F}_{L-}^{L}\right].$$
But now, from the optional stopping theorem,
$$\mathbb{E}\left[\widetilde{M}_{\infty
}\mid\mathcal{F}_{L}^{L}\right]=\widetilde{M}_{L};$$moreover, $M$
and $\widetilde{M}$ have the same jumps, and $L$ avoids
$\left(\mathcal{F}_{t}\right)$ stopping times, hence
$\widetilde{M}_{L}=\widetilde{M}_{L-},\;a.s.$ Hence,
$$\mathbb{E}\left[ \widetilde{M}_{\infty }\mid
\mathcal{F}_{L-}^{L}\right]=\widetilde{M}_{L}.$$Now, plugging this
into (\ref{e2}), and comparing with (\ref{e1}), we obtain the
equivalence between $(1)$ and $(2)$.
\end{proof}
\subsection{A solution related to Martingales which vanish at $L$}
It is a remarkable fact, discovered by Az\'{e}ma and Yor
(\cite{azemayorzero}), that a uniformly integrable martingale
vanishes at $L$, if and only if it is a solution to equation
(\ref{KMcond}):
\begin{prop}[Az\'{e}ma-Yor \cite{azemayorzero}]\label{orthogonalitparraportafl}
Let $M$ be an $L^{2}$ bounded martingale; then the following are
equivalent:
\begin{enumerate}
\item $\mathbb{E}\left[M_{\infty}\mid\mathcal{F}_{L}\right]=0$, or in
other words, $M_{\infty}\in
\left(L^{2}\left(\mathcal{F}_{L}\right)\right)^{\bot}$;
\item $M_{L}=0$.
\end{enumerate}
If one of the above conditions is satisfied, then $M\in
\mathcal{S}_{2}$.
\end{prop}
\begin{rem}
Proposition \ref{orthogonalitparraportafl} is still true if $M$ is
only assumed to be uniformly integrable.
\end{rem}
\begin{cor}\label{corilortho}
Let $M$ be an $L^{2}$ bounded martingale and let $$M^{L}\equiv
\mathbb{E}\left[M_{\infty}\mid\mathcal{F}_{L}\right].$$ Then the
martingale
$$M_{t}-\mathbb{E}\left[M^{L}\mid\mathcal{F}_{t}\right]$$belongs to
$\mathcal{S}_{2}$.
\end{cor}Theorem \ref{orthogonalitparraportafl} and the above corollary show
that it is enough to solve equation (\ref{KMcond}) when $M\in
L^{2}\left(\mathcal{F}_{L}\right)$. Indeed, $M_{\infty}$ can be
decomposed uniquely as: $$M_{\infty}=X_{1}+X_{2},$$where
$X_{1}\equiv M-M^{L} \in L^{2}\left(\mathcal{F}_{L}\right)^{\bot}$
and $X_{2}\equiv M^{L}\in L^{2}\left(\mathcal{F}_{L}\right)$. Thus,
$M_{t}=M_{t}^{1}+M_{t}^{2}$, with
$M_{t}^{1}=\mathbb{E}\left[X_{1}\mid\mathcal{F}_{t}\right]$,
$M_{t}^{2}=\mathbb{E}\left[X_{2}\mid\mathcal{F}_{t}\right]$, and
from Theorem \ref{orthogonalitparraportafl}, $M^{1}\in
\mathcal{S}_{2}$. Now, we give a description of $L^{2}$ martingales
$\left(M_{t}\right)$ such that $M_{\infty}=x_{L}$, where $(x_{t})$
is a predictable process (recall that we work under condition
$(\mathbf{A})$). Indeed, in all generality, for every $L^{2}$
martingale, there exists a predictable process $x$ such that:
$\mathbb{E}\left[M_{\infty}\mid\mathcal{F}_{L}\right]=x_{L}$.
\begin{prop}\label{propangulaire}
Let $(x_{t})$ be a predictable process such that
$\mathbb{E}\left[\left|x_{L}\right|\right]<\infty$. Then (for sake
of simplicity, we shall next write $A$ instead of $A^{L}$):
\begin{equation}\label{decoimport}
    \mathbb{E}\left[x_{L}\mid
    \mathcal{F}_{t}\right]=x_{L_{t}}\mathbb{P}\left(L\leq t\mid
    \mathcal{F}_{t}\right)+\mathbb{E}\left[\int_{t}^{\infty}x_{s}dA_{s}\mid
    \mathcal{F}_{t}\right],
\end{equation}where $$L_{t}=\sup\left\{s<t:\;1-Z_{s}^{L}=0\right\}.$$Moreover, the latter martingale can also be written as:
$$\mathbb{E}\left[x_{L}\mid
    \mathcal{F}_{t}\right]=-\int_{0}^{t}x_{L_{s}}d\mu_{s}^{L}+\mathbb{E}\left[\int_{0}^{\infty}x_{s}dA_{s}\mid
    \mathcal{F}_{t}\right],$$where $\left(\mu_{s}^{L}\right)$ is
    defined in Section 2 as the martingale part of the supermartingale $\left(Z_{t}^{L}\right)$.
\end{prop}
\begin{rem}
This proposition will be used in the next section to construct a
remarkable family of martingales.
\end{rem}
\begin{proof}
\begin{eqnarray*}
  \mathbb{E}\left[x_{L}\mid
    \mathcal{F}_{t}\right] &=& \mathbb{E}\left[x_{L}\mathbf{1}_{L\leq t}\mid
    \mathcal{F}_{t}\right]+\mathbb{E}\left[x_{L}\mathbf{1}_{L> t}\mid
    \mathcal{F}_{t}\right]  \\
   &=& x_{L_{t}}\mathbb{P}\left(L\leq t\mid
    \mathcal{F}_{t}\right)+\mathbb{E}\left[x_{L}\mathbf{1}_{L> t}\mid
    \mathcal{F}_{t}\right],
\end{eqnarray*}since from lemma \ref{azemalemmdeux}, on the set $\left\{L\leq t\right\}$, we have
$L_{t}=L$.
 Now, let $\Gamma_{t}$ be an
$\left(\mathcal{F}_{t}\right)$ measurable set;
$$\mathbb{E}\left[x_{L}\mathbf{1}_{L>
t}\mathbf{1}_{\Gamma_{t}}\right]=\mathbb{E}\left[\int_{t}^{\infty}x_{s}dA_{s}\mathbf{1}_{\Gamma_{t}}\right];$$hence
$$\mathbb{E}\left[x_{L}\mathbf{1}_{L> t}\mid
    \mathcal{F}_{t}\right]=\mathbb{E}\left[\int_{t}^{\infty}x_{s}dA_{s}\mid
    \mathcal{F}_{t}\right],$$and this completes the proof of the
    first part of the
    lemma. The second part follows from balayage arguments (see for example \cite{azemameyeryor}, Theorem 6.1); indeed:
    \begin{eqnarray*}
      x_{L_{t}}\mathbb{P}\left(L\leq t\mid
    \mathcal{F}_{t}\right) &=& x_{L_{t}}\left(1-Z_{t}^{L}\right) \\
       &=&
       -\int_{0}^{t}x_{L_{s}}d\mu_{s}^{L}+\int_{0}^{t}x_{s}dA_{s},
    \end{eqnarray*}where we have used the fact that $A$ lives on the set of times where $Z$ is equal to $1$. Now, since $\mathbb{E}\left[\int_{t}^{\infty}x_{s}dA_{s}\mid
    \mathcal{F}_{t}\right]=\mathbb{E}\left[\int_{0}^{\infty}x_{s}dA_{s}\mid
    \mathcal{F}_{t}\right]-\int_{0}^{t}x_{s}dA_{s}$, we have
    $$x_{L_{t}}\mathbb{P}\left(L\leq t\mid
    \mathcal{F}_{t}\right)+\mathbb{E}\left[\int_{t}^{\infty}x_{s}dA_{s}\mid
    \mathcal{F}_{t}\right]=-\int_{0}^{t}x_{L_{s}}d\mu_{s}^{L}+\mathbb{E}\left[\int_{0}^{\infty}x_{s}dA_{s}\mid
    \mathcal{F}_{t}\right],$$ and the proof of the lemma is now
    complete.
\end{proof}
Now, with the help of Proposition \ref{propangulaire}, we can solve
equation (\ref{KMcond}) for martingales of the form $M_{t}\equiv
\mathbb{E}\left[x_{L}\mid
    \mathcal{F}_{t}\right]$, and hence for any $L^{2}$ bounded
    martingale.
\begin{prop}\label{resolutiondeux}
Let $M_{t}\equiv \mathbb{E}\left[x_{L}\mid
    \mathcal{F}_{t}\right]$ be a uniformly integrable martingale ($(x_{t})$ is a predictable
    process). Then,
    \begin{eqnarray*}
      \mathbb{E}\left[M_{\infty}\mid
    \mathcal{F}_{L}\right] &=& \mathbb{E}\left[\int_{t}^{\infty}x_{s}dA_{s}\mid
    \mathcal{F}_{t}\right]|_{t=L} \\
       &=& \mathbb{E}\left[\int_{0}^{\infty}x_{s}dA_{s}\mid
    \mathcal{F}_{t}\right]|_{t=L}-\int_{0}^{\infty}x_{s}dA_{s}.
    \end{eqnarray*}Consequently, $\mathbb{E}\left[M_{\infty}\mid
    \mathcal{F}_{L}\right]=M_{L}$ if and only if $$\mathbb{E}\left[\int_{t}^{\infty}x_{s}dA_{s}\mid
    \mathcal{F}_{t}\right]|_{t=L}=x_{L}.$$
\end{prop}
\begin{rem}
We give in the next section some examples where all the calculations
can be done explicitly.
\end{rem}
\begin{proof}
This proposition is a consequence of Proposition \ref{propangulaire}
and the fact that: $A_{\infty}=A_{L}$.
\end{proof}

\subsection{Last zero before a fixed or a random time for the
standard Brownian Motion} We shall now use Proposition
\ref{orthogonalitparraportafl} to build martingales which are
solutions to equation (\ref{KMcond}) with a Brownian example which
has received much attention in the literature ( see \cite{zurich} or
\cite{azemajeulinknightyor} for more references). In the sequel, we
shall also use some results from \cite{Ashkanbessel}, where in
particular all the following results have been generalized to Bessel
processes of dimension $\delta(\equiv2(1-\mu))\in(0,2)$.

Let $\left( \Omega ,\mathcal{F},\left( \mathcal{F}_{t}\right) _{t\leq 1},%
\mathbf{P}\right) $ be a filtered probability space, where the filtration $%
\left( \mathcal{F}_{t}\right) $ is generated by a one dimensional
Brownian motion $\left( B_{t}\right) _{t\leq 1}$.
Let$$\gamma\equiv\sup\left\{t\leq1:\;B_{t}=0\right\}.$$

It is well known (see \cite{zurich, jeulinyor}) that:
\begin{equation}\label{zz1}
    \mathbb{P}\left[\gamma>t\mid
\mathcal{F}_{t}\right]=\sqrt{\dfrac{2}{\pi}}\int_{\frac{\left|B_{t}\right|}{\sqrt{1-t}}}^{\infty}\exp\left(\dfrac{-x^{2}}{2}\right)dx,\;t<1
\end{equation}
\begin{equation}\label{zz2}
    \lambda_{t}\equiv
A_{t}^{\gamma}=\sqrt{\dfrac{2}{\pi}}\int_{0}^{t}\dfrac{d\ell_{u}}{\sqrt{1-u}},
\;t<1.
\end{equation}
 Moreover, $m\equiv \dfrac{1}{\sqrt{1-\gamma}}\left|B_{1}\right|$,
$\varepsilon=sgn\left(B_{1}\right)$ and $\mathcal{F}_{\gamma}$ are
independent. We also have (as a consequence of Imhof's result, see
for example \cite{columbia} p.55):
$$\mathbb{P}\left(m\in d\rho\right)=\rho\exp\left(-\dfrac{\rho^{2}}{2}\right)d\rho.$$
\begin{rem}
A generalization of formulae (\ref{zz1}) and (\ref{zz2}) (which
leads to a multidimensional version of the arc sine law) is proved
in \cite{Ashkanbessel} for any Bessel process of dimension
$\delta\in (0,2)$.
\end{rem}
\begin{prop}[\cite{zurich}, chapter XIV]
Let $$M_{t}^{f}\equiv\mathbb{E}\left[f\left(B_{1}\right)\mid
\mathcal{F}_{t}\right]=P_{1-t}f\left(B_{t}\right),$$with $f$ a Borel
function such that:
$\mathbb{E}\left[\left|f\left(B_{1}\right)\right|\right]<\infty$,
and $\left(P_{t}\right)$ the semigroup of $\left(B_{t}\right)$. If
$f$ is an odd function, then $\left(M_{t}^{f}\right)$ is a solution
to equation (\ref{KMcond}), or in other words,
$$\mathbb{E}\left[M_{\gamma}^{f}\right]=\mathbb{E}\left[M_{\infty}^{f}\mid \mathcal{F}_{\gamma}\right].$$
\end{prop}
\begin{proof}
We have:
\begin{eqnarray}
  \mathbb{E}\left[f\left(B_{1}\right)\mid
\mathcal{F}_{\gamma}\right] &=& \frac{1}{2}\int_{-\infty}^{\infty}dx\left|x\right|\exp\left(-\dfrac{x^{2}}{2}\right)f\left(x\sqrt{1-\gamma}\right) \label{casbrown}\\
   &=&
   \dfrac{1}{2\left(1-\gamma\right)}\int_{0}^{\infty}dyy\exp\left(-\dfrac{y^{2}}{2\left(1-\gamma\right)}\right)\left(f\left(y\right)+f\left(-y\right)\right),
  \notag
\end{eqnarray}and hence, if $f$ is odd, then $\mathbb{E}\left[f\left(B_{1}\right)\mid
\mathcal{F}_{\gamma}\right]=0$, and from Proposition
\ref{orthogonalitparraportafl}, $M_{t}^{f}$ is a solution to
equation (\ref{KMcond}).
\end{proof}Now, let us consider the case when $f$ is an even
function such that
$\mathbb{E}\left[\left|f\left(B_{1}\right)\right|\right]<\infty$.
From (\ref{casbrown}), $$\mathbb{E}\left[f\left(B_{1}\right)\mid
\mathcal{F}_{\gamma}\right]=
\int_{0}^{\infty}dy\dfrac{y}{\left(1-\gamma\right)}\exp\left(-\dfrac{y^{2}}{2\left(1-\gamma\right)}\right)f\left(y\right).$$
Now, let us define:
$$M^{f,\gamma}\equiv f\left(B_{1}\right)-\int_{0}^{\infty}dy\dfrac{y}{\left(1-\gamma\right)}\exp\left(-\dfrac{y^{2}}{2\left(1-\gamma\right)}\right)f\left(y\right),$$
then it follows from Corollary \ref{corilortho} that the martingale
\begin{equation}\label{nouvellesmartes}
M_{t}^{f,\perp}\equiv\mathbb{E}\left[M^{f,\gamma}\mid
\mathcal{F}_{t}\right],\;t\leq1 \end{equation} is a solution to
equation (\ref{KMcond}). We first note that:
\begin{equation}\label{formulesemigroupr}
M_{t}^{f,\perp}=P_{1-t}f\left(B_{t}\right)-\int_{0}^{\infty}dyyf\left(y\right)\mathbb{E}\left[\dfrac{1}{\left(1-\gamma\right)}\exp\left(-\dfrac{y^{2}}{2\left(1-\gamma\right)}\right)\mid\mathcal{F}_{t}\right],
\end{equation}
and hence it is enough to have an explicit formula for martingales
of the form:
$$\mathbb{E}\left[h\left(\gamma\right)\mid\mathcal{F}_{t}\right],$$
where $h:[0,1]\rightarrow\mathbb{R}$ is a deterministic function.
This problem is solved in \cite{Ashkanbessel} (it suffices to take
$\mu=\frac{1}{2}$ to recover the Brownian setting) and leads to some
interesting results for our purpose.
\begin{lem}[\cite{Ashkanbessel}, with $\mu=\frac{1}{2}$]\label{decodemartarret}
Let $h:[0,1]\rightarrow\mathbb{R}_{+}$, be a Borel function, and let
$$\gamma\left(t\right)\equiv\sup\left\{u\leq t;\;B_{u}=0\right\}.$$Then:
\begin{equation*}
    \mathbb{E}\left[h\left(\gamma\right)\mid
    \mathcal{F}_{t}\right]=h\left(\gamma\left(t\right)\right)\left(1-Z_{t}^{\gamma}\right)+\mathbb{E}\left[h\left(\gamma\right)\mathbf{1}_{\left(\gamma>t\right)}\mid \mathcal{F}_{t}\right];
\end{equation*}with
\begin{equation}\label{equatioonpont}
    \mathbb{E}\left[h\left(\gamma\right)\mathbf{1}_{\left(\gamma>t\right)}\mid \mathcal{F}_{t}\right]=
   \dfrac{1}{\pi}\int_{0}^{1}dz\dfrac{h\left(t+z\left(1-t\right)\right)}{\sqrt{z\left(1-z\right)}}\exp\left(-\dfrac{B_{t}^{2}}{2z\left(1-t\right)}\right).
\end{equation}
\end{lem}
Now, with the help of Lemma \ref{decodemartarret}, after some
elementary calculations, we have the following  explicit expression
for the family of martingales $\left(M_{t}^{f,\perp}\right)$.
\begin{prop}
Let $f$ be an even Borel function such that
$\mathbb{E}\left[\left|f\left(B_{1}\right)\right|\right]<\infty$,
and let $$M_{t}^{f,\perp}=\mathbb{E}\left[
f\left(B_{1}\right)-\int_{0}^{\infty}dy\dfrac{y}{\left(1-\gamma\right)}\exp\left(-\dfrac{y^{2}}{2\left(1-\gamma\right)}\right)f\left(y\right)\mid\mathcal{F}_{t}\right],$$
as in (\ref{nouvellesmartes}). Define:
$$\theta\left(x\right)\equiv\sqrt{\dfrac{2}{\pi}}\int_{x}^{\infty}dv\exp\left(-\dfrac{v^{2}}{2}\right)=\int_{0}^{1}\dfrac{dv}{\pi\sqrt{v\left(1-v\right)}}\exp\left(-\dfrac{x^{2}}{2v}\right).$$
Then
$$M_{t}^{f,\perp}=M_{t}^{f,1}-M_{t}^{f,2}-M_{t}^{f,3},$$where:
\begin{eqnarray*}
  M_{t}^{f,1} &=& \int_{0}^{\infty}\dfrac{dz}{\sqrt{2\pi\left(1-t\right)}}f\left(z\right)\left(\exp\left(-\dfrac{\left(z+B_{t}\right)^{2}}{2\left(1-t\right)}\right)+\exp\left(-\dfrac{\left(z-B_{t}\right)^{2}}{2\left(1-t\right)}\right)\right), \\
  M_{t}^{f,2} &=& \theta\left(\dfrac{\left|B_{t}\right|}{\sqrt{1-t}}\right)\int_{0}^{\infty}dzzf\left(z\sqrt{1-\gamma_{t}}\right)\exp\left(-\dfrac{z^{2}}{2}\right) \\
  M_{t}^{f,3} &=&
  \int_{0}^{\infty}dzz\exp\left(-\dfrac{z^{2}}{2}\right)\int_{0}^{1}\dfrac{dw}{\pi\sqrt{w\left(1-w\right)}}f\left(z\sqrt{1-t}\sqrt{1-w}\right)\exp\left(-\dfrac{B_{t}^{2}}{2w\left(1-t\right)}\right),
\end{eqnarray*}and $\left(M_{t}^{f,\perp}\right)$ is a solution to
equation (\ref{KMcond}).
\end{prop}
\begin{rem}
The proposition shows that although Proposition
\ref{orthogonalitparraportafl} is elementary, it is in practice
difficult to compute the projection of the terminal value of a
martingale on the sigma algebra $\mathcal{F}_{\gamma}$.
\end{rem}
As a consequence of the explicit form for the martingales
$\mathbb{E}\left[h\left(\gamma\right)\mid\mathcal{F}_{t}\right]$, we
have the following first interesting result which shows how
(\ref{KMcond}) or (\ref{pdtacond}) may fail to hold in general:
\begin{prop}\label{corporgam}
Let $h:[0,1]\rightarrow\mathbb{R}_{+}$, be a Borel function, and
define $N_{t}^{h}=\mathbb{E}\left[h\left(\gamma\right)\mid
\mathcal{F}_{t}\right]$; then
\begin{equation}
    \mathbb{E}\left[N_{\infty}^{h}\mid
\mathcal{F}_{\gamma}\right]=h\left(\gamma\right),
\end{equation}whilst
\begin{equation}
    N_{\gamma}^{h}=\dfrac{1}{\pi}\int_{0}^{1}\dfrac{dv}{\sqrt{v\left(1-v\right)}}h\left(\gamma+v\left(1-\gamma\right)\right).
\end{equation}
\end{prop}\bigskip

The balayage formula can be used to get many solutions to equation
(\ref{KMcond}):
\begin{prop}
Define:
$$g_{t}\equiv\gamma\left(t\right)=\sup\left\{s\leq t:\;B_{s}=0\right\},$$and let $T>0$ be a
fixed time (thus with this notation, we have $g_{1}=\gamma$). Then,
for any bounded predictable process $\left(x_{s}\right)$,
$$X_{t}\equiv x_{g_{t\wedge T}}B_{t\wedge T}$$ is a uniformly
integrable martingale which satisfies (\ref{KMcond}) for $L=g_{T}$,
or more generally for $L=g_{t};\;t\leq T$.
\end{prop}
\begin{proof}
It is a consequence of the balayage formula that
$\left(x_{g_{t\wedge T}}B_{t\wedge T}\right)$ is a local martingale
(see \cite{revuzyor}, Chapter VI). From our assumptions, we easily
obtain that it is a bounded $L^{2}$ martingale. Now, $X_{g_{t}}=0$
for every $t\leq T$, and hence from Proposition
\ref{orthogonalitparraportafl}, $X$ satisfies (\ref{KMcond}).
\end{proof}
It is also possible to give many examples of honest times such that
the standard Brownian Motion, adequately stopped, satisfies
(\ref{KMcond}).
\begin{prop}
Let $T$ be a stopping time such that $\left(B_{t\wedge T}\right)$ is
a uniformly integrable martingale. Define $g_{T}$ as above. Then,
for every honest time $L\leq g_{T}$, we have $B_{L}=0$ and hence
$\left(B_{t\wedge T}\right)$ satisfies (\ref{KMcond}) for such
$L's$.
\end{prop}
\begin{proof}
As $L\leq g_{T}$, and both $L$ and $g_{T}$ are honest, we have
$\mathcal{F}_{L}\subset\mathcal{F}_{g_{T}}$. Consequently,
$$\mathbb{E}\left[B_{T}\mid\mathcal{F}_{L}\right]=\mathbb{E}\left[\mathbb{E}\left[B_{T}\mid\mathcal{F}_{g_{T}}\right]\mid\mathcal{F}_{L}\right]=0,$$
because $\mathbb{E}\left[B_{T}\mid\mathcal{F}_{g_{T}}\right]=0$ from
Proposition \ref{orthogonalitparraportafl}. Now, another application
of Proposition \ref{orthogonalitparraportafl} yields $B_{L}=0$ and
hence $\left(B_{t\wedge T}\right)$ satisfies (\ref{KMcond}) with
$L$.
\end{proof}
\begin{rem}
The last two propositions can be extended to continuous martingales.
\end{rem}
\section{Understanding the differences with a remarkable family of martingales}
So far, we have tried to characterize martingales for which, given a
random time, (\ref{pdtacond}) and (\ref{KMcond}) hold. Now, we try
to understand how these equalities may fail. Again, we consider the
case of honest times under condition $(\mathbf{A})$. To this end, we
introduce a family of uniformly integrable martingales, with some
remarkable properties, and which will serve us to test
(\ref{pdtacond}) and (\ref{KMcond}). From Corollary
\ref{relaorthogonalite} and Proposition
\ref{orthogonalitparraportafl}, it follows that interesting examples
of families of martingales such that (\ref{pdtacond}) and
(\ref{KMcond}) may fail to hold, should have the property:
$M_{\infty}$ is $\sigma\left(A_{\infty}\right)$ measurable (indeed,
$\sigma\left(A_{\infty}\right)\subset \mathcal{F}_{L}$ since
$A_{\infty}=A_{L}$). We should also mention that equation
(\ref{KMcond}) has been studied in the special case of David
Williams' pseudo-stopping time in \cite{AshkanYor}.
\subsection{A remarkable family of martingales}
We first prove a useful lemma:
\begin{lem}\label{leminermee}
Let $\varphi$ be a Borel function such that
$\mathbb{E}\left[\left|\varphi\left(A_{\infty}\right)\right|\right]<\infty$,
or equivalently
$\int_{0}^{\infty}dx\left|\varphi\left(x\right)\right|\exp\left(-x\right)<\infty$,
and let $\Phi\left(x\right)=\int_{0}^{x}dy\varphi\left(y\right)$.
Then:
\begin{equation}
\mathbb{E}\left[\varphi\left(A_{\infty}\right)\mid
\mathcal{F}_{t}\right]=\varphi\left(A_{t}\right)\left(1-Z_{t}\right)-\Phi\left(A_{t}\right)+\mathbb{E}\left[\Phi\left(A_{\infty}\right)\mid
\mathcal{F}_{t}\right].
\end{equation}
\end{lem}
\begin{proof}
From Lemma \ref{azemalemmdeux},
$\varphi\left(A_{\infty}\right)=\varphi\left(A_{L}\right)$, and
hence: $\mathbb{E}\left[\varphi\left(A_{\infty}\right)\mid
\mathcal{F}_{t}\right]=\mathbb{E}\left[\varphi\left(A_{L}\right)\mid
\mathcal{F}_{t}\right]$. We can thus apply Proposition
\ref{propangulaire} to obtain:
$$\mathbb{E}\left[\varphi\left(A_{\infty}\right)\mid
    \mathcal{F}_{t}\right]=\varphi\left(A_{L_{t}}\right)\mathbb{P}\left(L\leq t\mid
    \mathcal{F}_{t}\right)+\mathbb{E}\left[\int_{t}^{\infty}\varphi\left(A_{s}\right)dA_{s}\mid
    \mathcal{F}_{t}\right].$$ Now, from Lemma \ref{azemgeneral},
    $\varphi\left(A_{L_{t}}\right)=\varphi\left(A_{t}\right)$ and
    moreover, since $A$ is continuous,
    $$\int_{t}^{\infty}\varphi\left(A_{s}\right)dA_{s}=\int_{A_{t}}^{A_{\infty}}dx\varphi\left(x\right)=\Phi\left(A_{\infty}\right)-\Phi\left(A_{t}\right),$$and
    the assertion of the lemma follows.
\end{proof}
\begin{rem}\label{remreguliere}
If $f$ is a function of class $\mathcal{C}^{1}$, then, an
application of Lemma \ref{leminermee} with $\varphi=f'$ yields:
$$\mathbb{E}\left[f\left(A_{\infty}\right)-f'\left(A_{\infty}\right)\mid
\mathcal{F}_{t}\right]=f\left(A_{t}\right)-f'\left(A_{t}\right)\left(1-Z_{t}\right).$$
\end{rem}
Before introducing our family of martingales, we need to introduce
the following transform:  we associate with a continuous function
$\varphi$ the function $\widehat{\varphi}$, defined on
$\mathbb{R}_{+}$ by:
\begin{eqnarray*}
\widehat{\varphi}\left( x\right) &=&\exp \left( x\right)
\int_{x}^{\infty }dy\exp \left(
-y\right) \varphi\left( y\right) \\
&=&\int_{0}^{\infty }dy\exp \left( -y\right) \varphi\left(
y+x\right).
\end{eqnarray*}It is easy to see that $\widehat{\varphi}$ is a
function of class $\mathcal{C}^{1}$,
and:$$\widehat{\varphi}-\widehat{\varphi}'=\varphi.$$
\begin{prop}\label{mesmart}
Let $\varphi$ be a continuous function such that
$\mathbb{E}\left[\left|\varphi\left(A_{\infty}\right)\right|\right]<\infty$,
or equivalently
$\int_{0}^{\infty}dx\left|\varphi\left(x\right)\right|\exp\left(-x\right)<\infty$,
and let $$\widehat{\varphi}\left( x\right)=\exp \left( x\right)
\int_{x}^{\infty }dy\exp \left( -y\right) \varphi\left( y\right).$$
If
$\int_{0}^{\infty}dx\exp\left(-x\right)\left|\varphi\left(x\right)\right|x<\infty$,
then:
\begin{equation}\label{martremarkable}
\mathbb{E}\left[\varphi\left(A_{\infty}\right)\mid
\mathcal{F}_{t}\right]=Z_{t}\widehat{\varphi}\left(
A_{t}\right)+\left(1-Z_{t}\right)\varphi\left(A_{t}\right).
\end{equation}
\end{prop}
\begin{proof}
It suffices to apply Remark \ref{remreguliere} to
$\widehat{\varphi}$.
\end{proof}
\begin{rem}\label{remmonotone}
In fact, it can be shown, using monotone class arguments, that
formula \ref{martremarkable} remains valid if $\varphi$ is only
assumed to be a Borel function such that
$\int_{0}^{\infty}dx\left|\varphi\left(x\right)\right|\exp\left(-x\right)<\infty$.
These martingales have already been obtained by different means in
\cite{azemajeulinknightyor} and \cite{AshkanYorrem} (they are used
there in a different framework).
\end{rem}
One remarkable fact about these martingales, which we shall denote
by $\left(M_{t}^{\varphi}\right)$, is that we know their supremum
processes when $\varphi$\ is increasing. More precisely, we have:

\begin{prop}
\label{supremum}Assume that $\varphi$\ is a nonnegative, continuous
and increasing function such that the assumptions of Proposition
\ref{mesmart} are satisfied, and assume further that
$\left(M_{t}^{\varphi}\right)$ is a continuous martingale. Then the
supremum process of $\left(M_{t}^{\varphi}\equiv
\mathbb{E}\left[\varphi\left(A_{\infty}\right)\mid
\mathcal{F}_{t}\right]\right)$ is given by:%
\begin{equation*}
\sup_{s\leq t}M_{s}^{\varphi}=\widehat{\varphi}\left( A_{t}\right).
\end{equation*}
\end{prop}

\begin{proof}
We have:
\begin{equation}
\left( \widehat{\varphi}\left( A_{t}\right) -\varphi\left(
A_{t}\right) \right)\left(1-Z_{t}\right)
=-M_{t}^{\varphi}+\widehat{\varphi}\left( A_{t}\right) \label{abc}
\end{equation}%
Moreover,
\begin{equation*}
\widehat{\varphi}\left( A_{t}\right) -\varphi\left( A_{t}\right)
=\int_{0}^{\infty }dx\exp \left( -x\right) \left( \varphi\left(
x+A_{t}\right) -\varphi\left( A_{t}\right) \right)
\end{equation*}%
and since $\varphi$\ is increasing, we have $\left(
\widehat{\varphi}\left( A_{t}\right) -\varphi\left( A_{t}\right)
\right) \geq 0$. Similarly, we prove that $\widehat{\varphi}$\ is
increasing. Hence, (\ref{abc}) may be considered as a particular
case of Skorokhod's reflection equation
and thus:%
\begin{equation*}
\sup_{s\leq t}M_{s}^{\varphi}=\widehat{\varphi}\left( A_{t}\right)
\end{equation*}
\end{proof}\bigskip
\subsection{Martingales stopped at an honest time}
In the sequel, we assume that $\varphi$ is a continuous function
satisfying the following conditions of Proposition \ref{mesmart}:
$\int_{0}^{\infty}dx\left|\varphi\left(x\right)\right|\exp\left(-x\right)<\infty$,
and
$\int_{0}^{\infty}dx\exp\left(-x\right)\left|\varphi\left(x\right)\right|x<\infty$.
\begin{prop}
\label{maintm}Let $L$ be an honest time and
$M_{t}^{\varphi}=\mathbf{E}\left[ \varphi\left( A_{\infty
}\right) \mid \mathcal{F}_{t}\right] $. Then, we have:%
\begin{eqnarray}
M_{L}^{\varphi} &=&\exp \left( A_{L}\right) \int_{A_{L}}^{\infty
}dx\exp \left(
-x\right) \varphi\left( x\right)  \label{valeurterminale} \\
&=&\int_{0}^{\infty }dx\exp \left( -x\right) \varphi\left(
A_{L}+x\right) =\widehat{\varphi}\left( A_{L}\right)  \notag
\end{eqnarray}%
and
\begin{equation*}
\mathbf{E}\left[ \varphi\left( A_{\infty }\right) \mid
\mathcal{F}_{L}\right] =\varphi\left( A_{L}\right)
\end{equation*}
\end{prop}

\begin{proof}
The Proposition follows from Proposition \ref{mesmart}, and the fact
that $Z_{L}=1$, $A_{\infty}=A_{L}$.
\end{proof}

With Proposition \ref{maintm}, it is now clear why (\ref{KMcond})
may fail for honest times. More precisely;

\begin{cor}
\label{corimport}
\begin{equation*}
\mathbf{E}\left[ M_{\infty }^{\varphi}\mid \mathcal{F}_{L}\right]
=M_{L}^{\varphi}
\end{equation*}%
if and only if $\varphi$ is constant.
\end{cor}

\begin{proof}
If $\mathbf{E}\left[ M_{\infty }^{\varphi}\mid
\mathcal{F}_{L}\right] =M_{L}^{\varphi}$, then from Proposition
\ref{maintm}, $\varphi$ satisfies:
\begin{equation*}
\exp \left( y\right) \int_{y}^{\infty }dx\exp \left( -x\right)
\varphi\left( x\right) =\varphi\left( y\right) .
\end{equation*}%
The only solutions of this equation are the constant functions.
\end{proof}
\subsection{The expected value of martingales stopped at an honest
time} In the previous section, we saw that $\mathbf{E}\left[
M_{\infty }^{\varphi}\mid \mathcal{F}_{L}\right] $ and
$M_{L}^{\varphi}$ differ if the function $\varphi$ is not
constant. In this subsection, we shall compare the two quantities $\mathbf{E}%
\left[ M_{\infty }^{\varphi}\right] $ and $\mathbf{E}\left[
M_{L}^{\varphi}\right] $.
\begin{prop}
\label{expect}Let $L$\ be an honest time and
$M_{t}^{\varphi}=\mathbf{E}\left[ \varphi\left( A_{\infty
}\right) \mid \mathcal{F}_{t}\right]$. We have:%
\begin{eqnarray*}
\mathbf{E}\left[ M_{L}^{\varphi}\right] &=&\int_{0}^{\infty
}\int_{0}^{\infty
}dxdy\exp \left( -x\right) \exp \left( -y\right) \varphi\left( y+x\right) \\
&=&\mathbf{E}\left[ \varphi\left( \mathbf{e}_{1}+\mathbf{e}_{2}\right) \right] \\
&=&\int_{0}^{\infty }dx\exp \left( -x\right) x\varphi\left( x\right)
\end{eqnarray*}%
where $\mathbf{e}_{1}$ and $\mathbf{e}_{2}$\ are two independent
random variables following the standard exponential distribution,
whereas
\begin{eqnarray*}
\mathbf{E}\left[ M_{\infty }^{\varphi}\right] &=&\int_{0}^{\infty
}dx\exp \left(
-x\right) \varphi\left( x\right) \\
&=&\mathbf{E}\left( \varphi\left( \mathbf{e}_{1}\right) \right)
\end{eqnarray*}
\end{prop}

\begin{proof}
This is a consequence of Propositions \ref{maintm} and
\ref{mesmart}.
\end{proof}
\begin{prop}
If $\varphi$ is a positive increasing function, we have:%
\begin{equation*}
\sup_{t\geq 0}M_{t}^{\varphi}=M_{L}^{\varphi},
\end{equation*}and consequently, $\left( M_{t}^{\varphi}\right) \in
\mathcal{H}^{1}$ and $$\left\Vert M^{\varphi}\right\Vert
_{\mathcal{H}^{1}}=\mathbf{E}\left[ M_{L}^{\varphi}\right]. $$
\end{prop}

\begin{proof}
It is a consequence of Proposition \ref{supremum} and the fact that
$A_{\infty}=A_{L}$.
\end{proof}
But unlike the case of equation (\ref{KMcond}), we can find
solutions to equation (\ref{pdtacond}) among the martingales
$\left(M_{t}^{\varphi}\right)$. Recall that
$\left(\mathcal{L}_{n}\right)$ is the family of Laguerre
polynomials.
\begin{prop}
Let $L$\ be an honest time. Let $\varphi$ satisfy the conditions of
proposition \ref{mesmart}
 and $\mathbf{E}\left( \varphi^{2}\left( A_{\infty }%
\right) \right) <\infty $. Set again
$M_{t}^{\varphi}=\mathbf{E}\left[ \varphi\left( A_{\infty }\right)
\mid \mathcal{F}_{t}\right] $; then
\begin{equation*}
\mathbf{E}\left[ M_{\infty }^{\varphi}\right] =\mathbf{E}\left[
M_{L}^{\varphi}\right]
\end{equation*}%
\ if and only if
\begin{equation*}
\varphi\left( x\right) =\alpha+\sum_{n=2}^{\infty }\alpha
_{n}\mathcal{L}_{n}\left( x\right)
\end{equation*}%
where $\alpha\in \mathbb{R}$, and $\left( \alpha _{n}\right) $ are
such that $\sum \alpha _{n}^{2}<\infty $, i.e: in the development of
$\varphi$, the coefficients of $\mathcal{L}_{1} $ is $\alpha
_{1}=0$. In other words, $\varphi$\ belongs to the orthogonal of
$\mathcal{L}_{1}$.
\end{prop}

\begin{proof}
From Proposition \ref{expect}, $\mathbf{E}\left[ M_{\infty }^{\varphi}\right] =%
\mathbf{E}\left[ M_{L}^{\varphi}\right] $ if and only if:%
\begin{equation*}
\int_{0}^{\infty }dx\exp \left( -x\right) x\varphi\left( x\right)
=\int_{0}^{\infty }dx\exp \left( -x\right) \varphi\left( x\right)
\end{equation*}%
or equivalently:%
\begin{equation*}
\int_{0}^{\infty }dx\exp \left( -x\right) \mathcal{L}_{1}\left(
x\right) \varphi\left( x\right) =0
\end{equation*}%
It now suffices to develop $\varphi$\ in the basis $\left(
\mathcal{L}_{n}\left( x\right) \right) $ to conclude.
\end{proof}
\begin{proof}
It suffices to notice that $\left( \mathbf{E}\left[ M_{\infty }^{\varphi}\right] -%
\mathbf{E}\left[ M_{L}^{\varphi}\right] \right) $ is the coefficient
of $\mathcal{L}_{1}\left( x\right)=\left(
1-x\right) $\ in the expansion of $\varphi$\ in the basis $%
\left( \mathcal{L}_{n}\left( x\right) \right) $.
\end{proof}
\section*{acknowledgements}
I am deeply indebted to Marc Yor whose ideas and help were essential
for the development of this paper. I am also grateful to Thierry
Jeulin and Jan Obl\'{o}j for some fruitful discussions.

%\bibliographystyle{amsplain}
%\bibliography{theorgengros}

\begin{thebibliography}{99}
\bibitem{AARoy} \textsc{G.E. Andrews, R. Askey, R. Roy}: \textit{Special
functions}, Cambridge University Press (1999).

\bibitem{azema} \textsc{J. Az\'{e}ma}: \textit{Quelques applications de la th%
\'{e}orie g\'{e}n\'{e}rale des processus I}, Invent. Math.
\textbf{18} (1972) 293-336.

\bibitem{azemajeulinknightmokoyor} \textsc{J. Az\'ema, T. Jeulin, F. Knight, G. Mokobodzki, M. Yor}: \textit{Sur les
processus croissants de type injectif}, S\'{e}m.Proba. XXX, Lecture
Notes in Mathematics \textbf{1626}, (1996), 312-343.

\bibitem{azemajeulinknightyor} \textsc{J. Az\'ema, T. Jeulin, F. Knight, M. Yor}: \textit{Le th\'eor\`eme d'arr\^{e}t en une
fin d'ensemble pr\'evisible}, S\'{e}m.Proba. XXVII, Lecture Notes in
Mathematics \textbf{1557}, (1993), 133-158.

\bibitem{azemjeulknightyor} \textsc{J. Az\'{e}ma, T. Jeulin, F. Knight, M. Yor}: \textit{Quelques calculs de compensateurs
impliquant l'inj\'{e}ctivit\'{e} de certains processus croissants},
S\'{e}m.Proba. XXXII, Lecture Notes in Mathematics \textbf{1686},
(1998), 316-327.

\bibitem{azemameyeryor} \textsc{J. Az\'ema, P.A. Meyer, M. Yor}: \textit{Martingales relatives}, S\'{e}m.Proba. XXVI, Lecture Notes in Mathematics \textbf{1526},
(1992), 307-321.

\bibitem{azemayormartrem} \textsc{J. Az\'ema, M. Yor}: \textit{Etude d'une martingale remarquable}, S\'{e}m.Proba. XXIII, Lecture Notes in Mathematics \textbf{1372},
(1989), 88-130.

\bibitem{azemayorzero} \textsc{J. Az\'ema, M. Yor}: \textit{Sur les z\'{e}ros des martingales continues}, S\'{e}m.Proba. XXVI, Lecture Notes in Mathematics \textbf{1526},
(1992), 248-306.

\bibitem{barlow} \textsc{M.T. Barlow}, \textit{Study of a filtration
expanded to include an honest time}, ZW, \textbf{44}, 1978, 307-324.

\bibitem{beksy} \textsc{M.T. Barlow, M. Emery, F.B. Knight, S. Song, M. Yor}, \textit{Autour d'un th\'{e}or\`{e}me
de Tsirelson sur des filtrations browniennes et non browniennes},
S\'{e}m.Proba. XXXII, Lecture Notes in Mathematics \textbf{1686},
(1998), 264-305.

\bibitem{bremaudyor} \textsc{P. Br\'{e}maud, M. Yor}: \textit{Changes of
filtration and of probability measures}, Z.f.W, \textbf{45}, 269-295
(1978).

\bibitem{chungdoob} \textsc{K.L. Chung, J.L. Doob}, \textit{Fields, optionality and measurability},
Amer. J. Math, \textbf{87}, 1965, 397-424.

\bibitem{dellachmeyer} \textsc{C. Dellacherie, P.A. Meyer}: \textit{%
Probabilit\'{e}s et potentiel}, Hermann, Paris, vol II. 1980.

\bibitem{delmaismey} \textsc{C. Dellacherie, B. Maisonneuve, P.A. Meyer}:
\textit{Probabilit\'{e}s et potentiel}, Chapitres XVII-XXIV:
Processus de Markov (fin), Compl\'{e}ments de calcul stochastique,
Hermann (1992).

\bibitem{elliotjeanbyor} \textsc{R.J. Elliott, M. Jeanblanc, M. Yor}:
\textit{On models of default risk}, Math. Finance, \textbf{10},
179-196 (2000).

\bibitem{jeulin} \textsc{T. Jeulin}: \textit{Semi-martingales et
grossissements d'une filtration}, Lecture Notes in Mathematics
\textbf{833}, Springer (1980).

\bibitem{yorjeulin} \textsc{T. Jeulin, M. Yor}: \textit{Grossissement d'une
filtration et semimartingales: formules explicites}, S\'{e}m.Proba.
XII, Lecture Notes in Mathematics \textbf{649}, (1978), 78-97.

\bibitem{jeulinyor} \textsc{T. Jeulin, M. Yor (eds)}: \textit{Grossissements
de filtrations: exemples et applications}, Lecture Notes in
Mathematics \textbf{1118}, Springer (1985).

\bibitem{kazamaki} \textsc{N. Kazamaki}: \textit{Continuous exponential martingales and BMO},
Lecture Notes in Mathematics \textbf{1579}, Springer (1994).

\bibitem{Knight-Maisonneuve} \textsc{F.B. Knight, B. Maisonneuve}: \textit{A
characterization of stopping times}. Annals of probability,
\textbf{22}, (1994), 1600-1606.

\bibitem{Kunita-Watanabe} \textsc{H. Kunita, S. Watanabe}:
\textit{On square integarble martingales}. Nagoya J. Math.,
\textbf{30}, (1967), 209-245.

\bibitem{meyerdual} \textsc{P.A. Meyer}: \textit{Le dual de $H^{1}$ est $BMO$ (cas continu)}, S\'{e}m.Proba. VII,
Lecture Notes in Mathematics \textbf{321}, (1973), 136-145.

\bibitem{columbia} \textsc{R. Mansuy, M. Yor}: \textit{Random times and (enlargement of filtrations) in a Brownian setting},
Lecture Notes in Mathematics, \textbf{1873}, Springer (2006).

\bibitem{AshkanYor} \textsc{A. Nikeghbali, M. Yor}: \textit{A definition and
some characteristic properties of pseudo-stopping times},  Ann.
Prob. \textbf{33}, (2005) 1804-1824.

\bibitem{ashyordoob} \textsc{A. Nikeghbali, M. Yor}: \textit{Doob's maximal identity, multiplicative decompositions
 and enlargements of filtrations}, to appear in  the Illinois Journal of Mathematics and available on the ArXiv.

 \bibitem{Ashkanbessel} \textsc{A. Nikeghbali}: \textit{Some random times and martingales associated with
$BES_{0}(\delta)$ processes $(0<\delta<2)$}, Alea \textit{2} 67-89
(2006).

\bibitem{AshkanYorrem} \textsc{A. Nikeghbali, M. Yor}: \textit{A class of remarkable submartingales}, Stochastic Processes and their Applications \textbf{116}, (2006) 917-938.

\bibitem{Ashkanassay} \textsc{A. Nikeghbali}: \textit{An essay on the general theory of stochastic processes}, to appear in Probability Surveys, preliminary version available on the ArXiv.

\bibitem{pitmanyor} \textsc{J.W. Pitman, M. Yor}: \textit{Bessel processes
and infinitely divisible laws, In: D. Williams (ed.) Stochastic
integrals}, Lecture Notes in Mathematics \textbf{851}, Springer
(1981).

\bibitem{protter} \textsc{P.E. Protter}: \textit{Stochastic integration and
differential equations}, Springer. Second edition, version 2.1
(2005).

\bibitem{revuzyor} \textsc{D. Revuz, M. Yor}: \textit{Continuous martingales
and Brownian motion}, Springer. Third edition (1999).

\bibitem{rogerswilliams} \textsc{C. Rogers, D. Williams}: \textit{%
Diffusions, Markov processes and Martingales, vol 2: Ito calculus},
Wiley and Sons, New York, 1987.

\bibitem{williams} \textsc{D. Williams}: \textit{A non stopping time with
the optional-stopping property}, Bull. London Math. Soc. \textbf{34}
(2002), 610-612.

\bibitem{gros} \textsc{M. Yor}: \textit{Grossissement d'une filtration et semi-martingales: th\'eor\`emes g\'en\'eraux}, S\'{e}m. Proba. XII, Lecture Notes in Mathematics \textbf{649},
(1978).

\bibitem{zurich} \textsc{M. Yor}: \textit{Some aspects of Brownian motion},
\textit{Part II. Some recent martingale problems.} Birkhauser, Basel
(1997).
\end{thebibliography}
\end{document}